\documentclass[english]{article}

\usepackage{subfigure}
\usepackage{graphicx}

\usepackage[T1]{fontenc}
\usepackage[latin9]{inputenc}
\usepackage{float}
\usepackage{bm}

\usepackage{amsmath}
\usepackage{amsthm,amssymb,amsfonts}
\usepackage{upgreek}
\usepackage[hmargin=1.2in]{geometry}
\usepackage[textsize=footnotesize]{todonotes}

\makeatletter
\providecommand{\tabularnewline}{\\}
\makeatother
\usepackage{babel}

\begin{document}

\newcommand{\bx}{\mathbf{x}}
\newcommand{\bu}{\mathbf{u}}
\newcommand{\bq}{\mathbf{q}}
\newcommand{\bqh}{\mathbf{\hat{q}}}
\newcommand{\bg}{\mathbf{g}}
\newcommand{\bff}{\mathbf{f}}
\newcommand{\br}{\mathbf{r}}
\newcommand{\mA}{\mathbf{A}}
\newcommand{\mB}{\mathbf{B}}
\newcommand{\bQ}{\mathbf{Q}}
\newcommand{\bF}{\mathbf{F}}
\newcommand{\bG}{\mathbf{G}}
\newcommand{\bZero}{\mathbf{0}}

\newcommand{\DQ}{\Delta \mathbf{Q}}
\newcommand{\Dq}{\Delta \mathbf{q}}
\newcommand{\Dx}{\Delta x}
\newcommand{\Dy}{\Delta y}

\newcommand{\Zwave}{\mathbf{\cal Z}}
\newcommand{\amdQ}{{\cal A^-}\DQ}
\newcommand{\apdQ}{{\cal A^+}\DQ}
\newcommand{\bmdQ}{{\cal B^-}\DQ}
\newcommand{\bpdQ}{{\cal B^+}\DQ}
\newcommand{\amdq}{{\cal A^-}\Dq}
\newcommand{\apdq}{{\cal A^+}\Dq}
\newcommand{\bmdq}{{\cal B^-}\Dq}
\newcommand{\bpdq}{{\cal B^+}\Dq}
\newcommand{\adq}{{\cal A}\Dq}
\newcommand{\bdq}{{\cal B}\Dq}

\newcommand{\KA}{K_\textup{A}}
\newcommand{\KB}{K_\textup{B}}
\newcommand{\rhoA}{\rho_\textup{A}}
\newcommand{\rhoB}{\rho_\textup{B}}
\newcommand{\half}{\frac{1}{2}}
\newcommand{\imh}{{i-\half}}
\newcommand{\iph}{{i+\half}}
\newcommand{\jmh}{{j-\half}}
\newcommand{\jph}{{j+\half}}

\title{Numerical simulation of cylindrical solitary waves in periodic media}

\author{
  Manuel Quezada de Luna
  \thanks{Texas A\&M University, College Station, TX 77843
  (\mbox{mquezada@math.tamu.edu})}
   \and
   David I. Ketcheson
   \thanks{King Abdullah University of Science and Technology,
    Box 4700, Thuwal, Saudi Arabia, 23955-6900
    (\mbox{david.ketcheson@kaust.edu.sa}) }
  }
  
\maketitle

\abstract{
        We study the behavior of nonlinear waves in a two-dimensional
        medium with density and stress relation that vary periodically in space.  
        Efficient approximate Riemann solvers are developed for the corresponding
        variable-coefficient first-order hyperbolic system.  We present direct
        numerical simulations of this multiscale problem, focused on the
        propagation of a single localized perturbation
        in media with strongly varying impedance.
        For the conditions studied, we find little evidence of shock formation.  Instead, 
        solutions consist primarily of solitary waves.
        These solitary waves are observed
        to be stable over long times and to interact in a manner approximately
        like solitons.
        The system considered has no dispersive terms; these
        solitary waves arise due to the material heterogeneity, which leads to 
        strong reflections and effective dispersion.
        }

\section{Introduction}\label{sec:introduction}

Solutions of first order hyperbolic partial
differential equations (PDEs) generically
lead to formation of shock singularities and subsequent entropy decay.
In contrast,
nonlinear wave equations with dispersive terms, such as the Korteweg-deVries
equation (KdV), may exhibit solitary wave 
solutions \cite{zabusky1965interaction} . These solutions are remarkable in the
context of nonlinear PDEs since
they propagate without changing shape and, in many cases, interact only through
a phase shift.  Thus the long-time solution behavior of nonlinear waves depends
critically on the presence of dispersive regularizing terms.

Santosa and Symes \cite{santosa1991dispersive} showed that solutions of the
linear wave equation in a medium with periodically varying coefficients
exhibit dispersive behavior, even though the equation contains no
dispersive terms.  LeVeque and Yong \cite{leveque2003solitary} found that
solutions of the 1D $p$-system in a periodic layered
medium form solitary waves (referred to as {\em stegotons} due to their discontinuous
shape) similar to those arising in
nonlinear dispersive wave equations like KdV.  This is remarkable since the
equations they considered are first order hyperbolic PDEs
with no dispersive terms.  Recently it has been observed that solutions of
first order hyperbolic systems with periodic coefficients may be characterized
by shocks or by solitary waves, depending on the medium and the initial
conditions \cite{Ketcheson_LeVeque_2011}.
Further experiments in \cite{ketchesonphdthesis,Ketcheson_LeVeque_2011} indicate
that solitary waves may arise quite generally in the solution of nonlinear hyperbolic
PDEs with periodically varying coefficients.
Those results have revealed new behaviors of nonlinear waves in 
1D heterogeneous materials, and it is natural to ask if such behaviors
persist in higher dimensions.  The present work, which builds on
\cite{Quezada_de_Luna_2011}, extends the aforementioned studies to a fully two-dimensional setting. 

We consider nonlinear wave propagation in a 2D periodic medium, modeled by
variable-coefficient first order hyperbolic PDEs.  Since 1D solitary
waves are often found to be unstable when extended to higher dimensions,
a principal question is whether waves like stegotons are stable in 2D, and
whether they form in materials that vary in multiple spatial directions.
In order to investigate these questions in a general setting, 
we have chosen to study the nonlinear wave equation
\begin{align*}
  \epsilon_{tt}-\nabla\cdot\left(\frac{1}{\rho(\bx)} \nabla\sigma(\epsilon,\bx)\right) & = 0.
\end{align*}
This is perhaps the simplest multidimensional nonlinear wave model that is general in 
the sense of allowing for wave propagation in all directions.
In 1D, it is commonly referred to as the $p$-system due to its connection
with Lagrangian gas dynamics.

\begin{figure}
 \begin{centering}
  \subfigure[Initial condition (close-up)]{\includegraphics[width=2.5in]{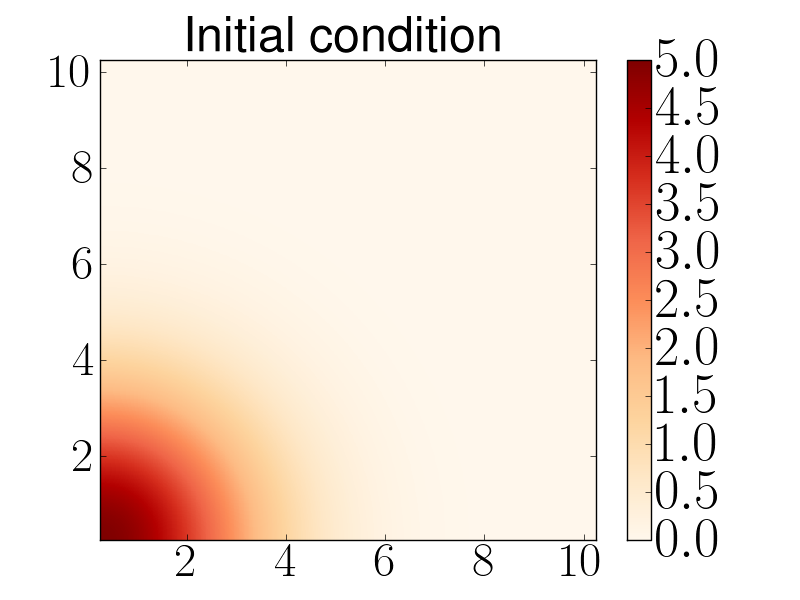}\label{fig:IC}}
  \subfigure[Checkerboard medium]{\includegraphics[scale=0.3]{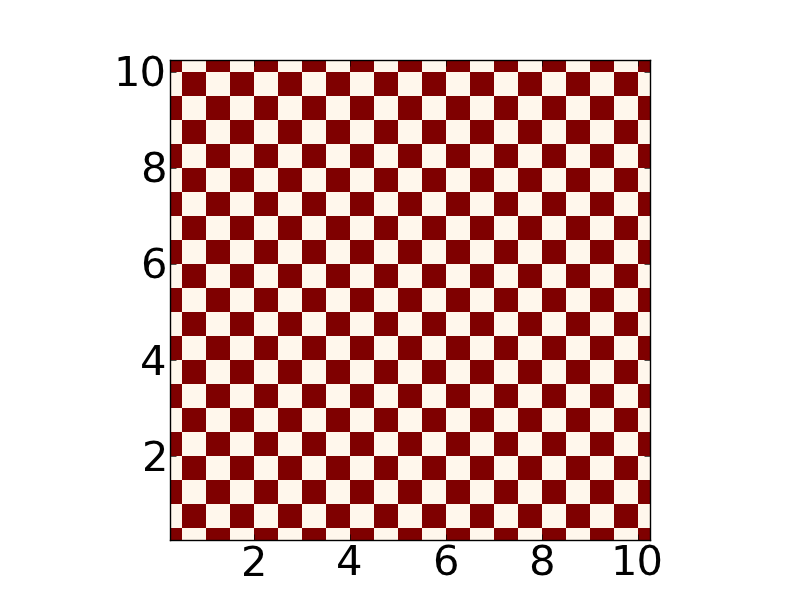}\label{Flo: checker domain1}}
  \subfigure[Solution in four different media]{\includegraphics[width=4in]{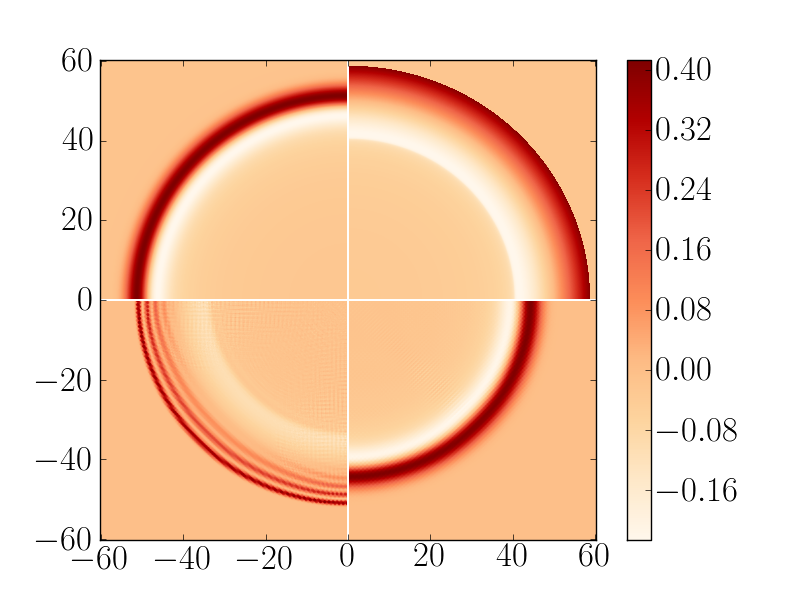}\label{fig:quadrants}}
  \caption{Wave behavior in different media.  Top left: linear, homogeneous.  
           Top right: nonlinear, homogeneous.  Bottom right: linear, periodic.
           Bottom left: nonlinear, periodic. The quantity plotted is stress ($\sigma$).
           }
 \end{centering}
\end{figure}

In order to introduce the object of our study, four representative examples of
the behavior of cylindrical wavefronts are depicted in Figure \ref{fig:quadrants}.  
Each solution shown results from the same initial condition: a small,
cylindrically symmetric Gaussian pulse centered at the origin,
shown (close-up) in Figure \ref{fig:IC} and given by: 
\begin{equation} \label{eq: initial condition}
  \sigma(\bm{x},0)=5\exp\left(-\frac{(x-0.25)^2}{10}-\frac{(y-0.25)^2}{10}\right).
\end{equation}
Since the solution is symmetric under reflection about the $x$- and $y$-axes,
a single quadrant is sufficient to characterize the full solution.
The top left quadrant of Figure \ref{fig:quadrants} shows the solution
obtained in a linear, homogeneous medium:
a smooth pulse expanding at the sound speed of the medium.  The top right
quadrant corresponds to a nonlinear, homogeneous medium; the leading
edge of the pulse steepens into a shock wave.  The bottom right quadrant
corresponds to a linear, heterogeneous medium shown in Figure \ref{Flo: checker domain1} 
that is composed of alternating
square homogeneous regions in a checkerboard pattern.
The front travels more slowly than in the homogeneous
case, due to the effect of reflections.  Finally, the bottom left quadrant
shows the subject of the present study.  The medium is both nonlinear and
heterogeneous (with the same checkerboard structure).  In this case,
a train of cylindrical pulses forms. As we will see, these pulses appear to behave
as solitary waves and we refer to them as cylindrical stegotons.
Note that they do not exhibit exact cylindrical symmetry, due to the lack
of such symmetry in the medium.

In the remainder of the paper, we investigate the conditions that lead to these
{\em cylindrical stegotons}, as well as their properties and behavior.  
The model equations and materials are introduced in Section \ref{sec:psys}.

Accurate numerical modeling of nonlinear waves in multidimensional heterogeneous
media is challenging.  Care must be taken to properly handle material
inhomogeneities (especially discontinuities) and shocks.
We use finite volume methods based on approximate Riemann solvers
implemented in Clawpack
\cite{LeVeque_Berger_2011} and SharpClaw \cite{Ketcheson_Parsani_LeVeque_2011},
New approximate Riemann solvers for the heterogeneous 2D $p$-system
are discussed in Section \ref{sec:methods}.
Because very large grids (with approximately 7 billion unknowns) 
are required to resolve this multiscale problem over times
and distances sufficient to observe solitary wave formation,
we have used the parallel PyClaw framework 
\cite{pyclaw-sisc}
in order to run on 16,384 cores of the Shaheen system at KAUST.

Motivated by known
behaviors of more traditional solitary waves, and also by studies of 1D
stegotons, we investigate qualitatively the following questions:
\begin{itemize}
    \item What kinds of media and initial conditions give rise to solitary waves?
    \item How do cylindrical stegotons interact with each other?
    \item What is the role of shock formation in 2D periodic media with strongly
            varying impedance?
\end{itemize}
These questions are addressed through further numerical experiments in Section
\ref{sec:results}.

\section{The 2D spatially-varying $p$-system\label{sec:psys}}

    \subsection{Governing equations}

As mentioned in the introduction, we have chosen to study the $p$-system due to
its relative simplicity and generality.  This system can be viewed as a
simplification of models governing more complex wave behavior, such as that of
elastic waves. In a solid, an elastic wave is composed of longitudinal or
P-waves and transversal or S-waves \cite{leveque2002finite}. If we assume
the stress is hydrostatic (i.e., there is no shear stress and the extensional
stress components are equal), the propagation of elastic waves can be modeled by:
\begin{align}\label{psystem-2nd-order}
  \epsilon_{tt}-\nabla\cdot\left(\frac{1}{\rho(\bx)} \nabla\sigma(\epsilon,\bx)\right) & = 0,
\end{align}
where $\epsilon$ represents the strain, $\rho(\bx)$ is the spatially-varying
material density, $\sigma(\epsilon,\bx)$ is the stress and $\bx=[x,y]^{T}$ is the position vector. 
Similar to \cite{leveque2003solitary}, we consider
the nonlinear constitutive relation
  \begin{align}
    \sigma(\epsilon,\bx) & = \exp(K(\bx)\epsilon)+1, \label{eq: nonlinear const relation}
  \end{align}
where $K(\bx)$ plays the role of bulk modulus. 
Equation \eqref{psystem-2nd-order} with the stress relation \eqref{eq: nonlinear const relation}
admits shock formation.  In order to determine entropy-satisfying weak
solutions, let us write \eqref{psystem-2nd-order} as a first-order
hyperbolic system of conservation laws:
\begin{subequations} \label{psystem}
  \begin{equation}
    \bq_{t}+\bff(\bq,\bx)_x+\bg(\bq,\bx)_{y}=\bZero,\label{eq: p system conservation form}\end{equation}
where
  \begin{align}
    \bq & = \begin{bmatrix} \epsilon\\ \rho(\bx) u \\ \rho(\bx)  v \end{bmatrix},
    & \bff(\bq,\bx) & =\begin{bmatrix} u \\ -\sigma(\epsilon,\bx)\\ 0\end{bmatrix},
    & \bg(\bq,\bx)=\begin{bmatrix} v\\ 0\\ -\sigma(\epsilon,\bx)\end{bmatrix}.
  \end{align}
\end{subequations}
Here $u$ and $v$ are the $x$- and $y$-components of velocity,
$\bq$ is the vector of conserved quantities, and $\bff,\bg$ 
are the components of the flux in the $x$- and $y$-directions, respectively. 
This form arises naturally through consideration
of kinematics and Newton's second law, and leads to the correct jump conditions
across shocks.

  \subsection{Periodic Media \label{sub: Medium}}
\begin{figure}
 \begin{centering}
  \subfigure[Checkerboard medium (equation \eqref{checkerboard})]{\includegraphics[scale=0.3]{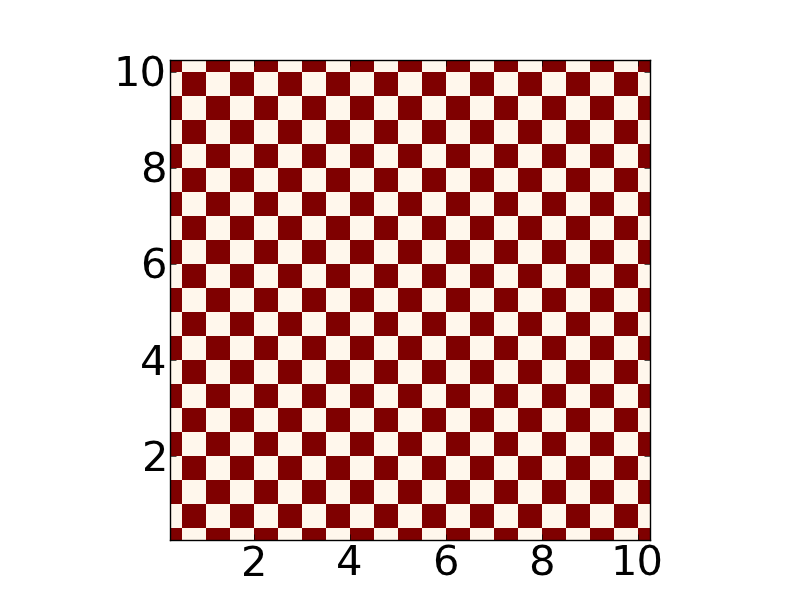}\label{Flo: checker domain}}
  \subfigure[Sinusoidal type medium (equation \eqref{sinusoidal})]{\includegraphics[scale=0.3]{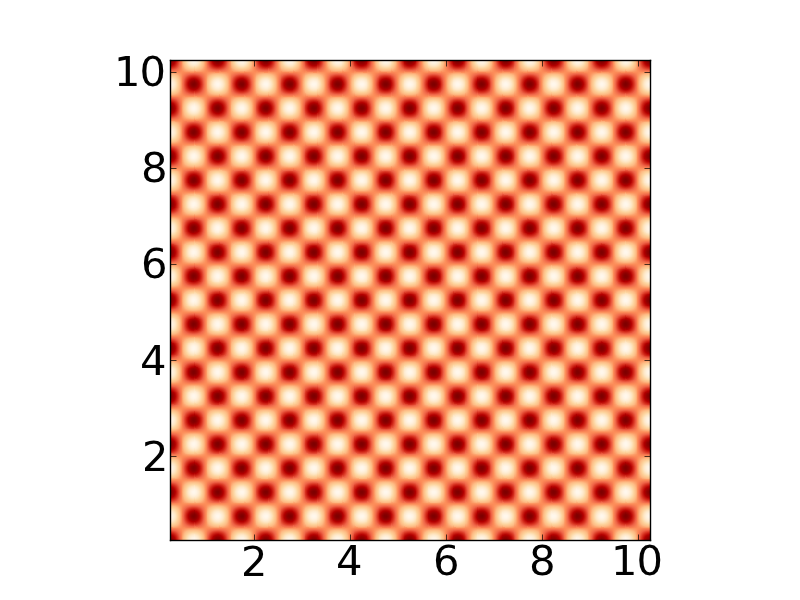}\label{Flo: sinusoidal domain}}
 \caption{The domains discussed in Section \ref{sub: Medium}.  Sections of size $10\times 10$ are 
 shown here, but much larger domains are used in the numerical simulations.\label{Flo: domains}}
 \end{centering}
\end{figure}

We consider two types of periodic variation in the density
$\rho(\bx)$ and the bulk modulus $K(\bx)$.  The period of the medium is always
taken to be unity.
The first medium, shown in figure \ref{Flo: checker domain}, consists of a
checkerboard pattern:
\begin{align} \label{checkerboard}
  K(\bx),\rho(\bx) & = \begin{cases}
                    (K_A,\rho_A) \mbox{ if } \left(x-\lfloor x\rfloor-\half\right)
                    \left(y-\lfloor y\rfloor-\half\right)<0 \\
                    (K_B,\rho_B) \mbox{ if } \left(x-\lfloor x\rfloor-\half\right)
                    \left(y-\lfloor y\rfloor-\half\right)>0
  \end{cases}
\end{align}
The second medium, shown in figure \ref{Flo: sinusoidal domain}, is similar
but smoothly (sinusoidally) varying:
\begin{subequations} \label{sinusoidal}
\begin{align}
  K(\bx) & = \frac{\KA+\KB}{2} +
                \frac{\KA-\KB}{2} \sin\left(2\pi x\right)\sin\left(2\pi y\right) \\
  \rho(\bx) & = \frac{\rhoA+\rhoB}{2} +
                \frac{\rhoA-\rhoB}{2} \sin\left(2\pi x\right)\sin\left(2\pi y\right)
\end{align}
\end{subequations}
We always take $K_A=\rho_A=1$.

\section{Numerical Discretization\label{sec:methods}}
In this section we describe the numerical methods used to solve the 2D
spatially varying $p$-system \eqref{psystem}.  For all
computations in this work we use the algorithms implemented in Clawpack 
and SharpClaw.  Both are finite volume methods for solving hyperbolic PDEs
based on solving Riemann problems.  
For the simulations presented here, both Clawpack and SharpClaw are
driven through PyClaw \cite{pyclaw-sisc}
which is a lightweight Python
framework that calls the low-level Fortran routines of Clawpack
and SharpClaw, and also interfaces with PETSc \cite{petsc} to provide
parallelism.

After briefly introducing the discretization schemes, we focus on their
key ingredient, the approximate solution of the Riemann problem for the
spatially-varying $p$-system \eqref{psystem}.
The Riemann solvers developed here follow the ideas of
\cite{leveque2002finitePaper}.  In particular, they are based on
use of an all-shock solution and $f$-wave decomposition \cite{bale2003wave}.

\subsection{Second-order Clawpack discretization}
Clawpack is based on Lax-Wendroff discretization combined with TVD limiters,
and is second order accurate in space and time \cite{LeVeque_Berger_2011}.  
The multidimensional Clawpack algorithms used here require the propagation of waves in both
the normal and transverse directions at each cell edge.
The Clawpack discretization takes the form 
\begin{align*}
  \bQ_{ij} & = \bQ_{ij}^{n}-\frac{\Delta t}{\Delta x}\left({\cal A}^{+}\Delta\bQ_{i-\frac{1}{2},j}^{n}+{\cal A}^{-}\Delta\bQ_{i+\frac{1}{2},j}^{n}\right)\\
   & -\frac{\Delta t}{\Delta y}\left({\cal B}^{+}\Delta\bQ_{i,j-\frac{1}{2}}^{n}+{\cal B}^{-}\Delta\bQ_{i,j+\frac{1}{2}}^{n}\right)\\
  &  -\frac{\Delta t}{\Delta x}\left(\tilde{\bm{F}}_{i+\frac{1}{2},j}^{n}-\tilde{\bm{F}}_{i-\frac{1}{2},j}^{n}\right)-\frac{\Delta t}{\Delta y}\left(\tilde{\bm{G}}_{i,j+\frac{1}{2}}^{n}-\tilde{\bm{G}}_{i,j-\frac{1}{2}}^{n}\right),
\end{align*}
where ${\cal A}^{\pm}\Delta\bQ_{i\mp\frac{1}{2},j}^{n}$ and 
${\cal B}^{\pm}\Delta\bQ_{i,j\mp\frac{1}{2}}^{n}$ are first-order
fluctuations computed by the normal Riemann solvers described in 
Section \ref{sec:normal}.
The quantities $\tilde{\bm{F}}_{i\pm\frac{1}{2},j}^{n}$
and $\tilde{\bm{G}}_{i,j-\frac{1}{2}}^{n}$ are second-order corrections that
include the fluctuations computed by the transverse Riemann solvers
described in Section \ref{sec:transverse}, as well as
high-resolution approximations to $\bq_{xx},\bq_{yy}$.
For details, see \cite{leveque2002finite}. 

\subsection{High-order WENO discretization using SharpClaw}
SharpClaw is based on the method of lines approach and uses fifth-order
weighted essentially non-oscillatory (WENO) reconstruction in space with fourth-order
strong stability preserving (SSP) Runge-Kutta integration in time 
\cite{Ketcheson_Parsani_LeVeque_2011}.
SharpClaw requires propagation of waves only in the direction normal to
each edge.

    First, a WENO reconstruction of the solution is computed
    from the cell averages $Q_i$ to give high order accurate point values
    $q^\textup{L}_\imh$ and $q^\textup{R}_\imh$ just to the left and right
    (respectively) of each cell interfaces $x_\imh$.  
    A Riemann solution is computed at each interface based on the reconstructed
    values there.
    The resulting fluctuations are used to update the adjacent cell averages.
    An additional term appears that is proportional to $\int_{x_\imh}^{x_\iph}Aq_x dx$.
    For conservative systems like \eqref{psystem}, this term can be conveniently computed in terms
    of a fictitious internal Riemann problem in each cell \cite{Ketcheson_Parsani_LeVeque_2011}.
    The semi-discrete scheme takes the form
    \begin{align*}
    \frac{\partial \bQ^{n+1}_{ij}}{\partial t} = 
            - \frac{1}{\Dx}\left(\amdq_{\iph,j} + \apdq_{\imh,j} + \adq_{ij}\right)
            - \frac{1}{\Dy}\left(\bmdq_{i,\jph} + \bpdq_{i,\jmh} + \bdq_{ij}\right)
    \end{align*}

    SharpClaw employs the same wave propagation Riemann solvers and user
    interface as Clawpack.
    In multi-dimensions SharpClaw requires propagation of waves only in the
    normal direction to each edge.

\subsection{Normal Riemann solvers\label{sec:normal}}
We assume that the density and bulk modulus are constant
within each grid cell: $\rho=\rho_{ij}, K=K_{ij}$. 
In the case of smoothly varying media, this is an approximation.
Then the system of conservation laws \eqref{psystem} 
can be written in the following quasilinear form within cell $i,j$:
\begin{align}
\bq_t + \mA_{ij} \bq_x + \mB_{ij} \bq_y & = \bZero
& \forall (x,y) \in (x_\imh,x_\iph)\times(y_\jmh,y_\jph)
\end{align}
where
\begin{align}
\bq & = \begin{bmatrix} \epsilon \\ \rho u \\ \rho v
        \end{bmatrix}, &
\mA_{ij} = \bff'_{ij}(\bq) & = \begin{bmatrix}
            0 & -\frac{1}{\rho_{ij}} & 0 \\
            -\sigma_{\epsilon,ij} & 0 & 0 \\
            0 & 0 & 0
        \end{bmatrix}, &
\mB_{ij} = \bg'_{ij}(\bq) & = \begin{bmatrix}
            0 & 0 & -\frac{1}{\rho_{ij}} \\
            0 & 0 & 0 \\
            -\sigma_{\epsilon,ij} & 0 & 0
        \end{bmatrix}. &
\end{align}
Note that $\sigma_{\epsilon,ij}$ denotes the derivative of $\sigma(\epsilon,x_i,y_j)$
with respect to $\epsilon$.
For concreteness, we consider a Riemann problem in the x-direction,
at interface $x_{\imh,j}$ which consists of the hyperbolic system
\eqref{psystem} with coefficients
\begin{align}
 \rho(x,y) & = \begin{cases} \rho_{i-1,j} & \forall x<x_\imh \\ 
                           \rho_{i,j} & \forall x>x_\imh \end{cases}, &
 K(x,y) & = \begin{cases} K_{i-1,j} & \forall x<x_\imh \\ 
                           K_{i,j} & \forall x>x_\imh; \end{cases}
\end{align}
and initial condition
\begin{align}
 \bq(x,y) & = \begin{cases} \bQ_{i-1,j} & \forall x<x_\imh \\ 
                           \bQ_{i,j} & \forall x>x_\imh \end{cases}.
\end{align}
In practice, $\bQ_{ij}$ may represent a cell-average (in Clawpack) or a reconstructed
value at the cell interface (in SharpClaw).
The eigenvectors of $\mA_{ij}(\bQ_{ij})$ are
\begin{align*}
  \br_{ij}^{1} & = \begin{bmatrix} 1 \\ Z_{ij} \\ 0 \end{bmatrix}, & 
  \br_{ij}^{2} & = \begin{bmatrix} 0 \\ 0       \\ 1 \end{bmatrix}, &
  \br_{ij}^{3} & = \begin{bmatrix}-1 \\ Z_{ij} \\ 0 \end{bmatrix}
\end{align*}
with corresponding eigenvalues $\{-c_{ij},0,+c_{ij}\}$.  Here
$Z_{ij}=\sqrt{\rho_{ij}\sigma_{\epsilon,ij}}$ is the impedance and 
$c_{ij}=\sqrt{\frac{\sigma_{\epsilon,ij}}{\rho_{ij}}}$ is the sound speed.

In the linear case, each wave in the Riemann solution is a discontinuity
proportional to the corresponding eigenvector {\em in the material carrying
the wave}.  Thus the solution can be found by decomposing the
difference $\DQ_{\imh,j} = \bQ_{ij}-\bQ_{i-1,j}$ in terms of the following three eigenvectors:
\begin{align} \label{imh_eigenvectors}
  \br_{\imh,j}^1 = \br_{i-1,j}^{1} & = \begin{bmatrix} 1 \\ Z_{i-1,j} \\ 0 \end{bmatrix}, & 
  \br_{\imh,j}^2 = \br^{2}         & = \begin{bmatrix} 0 \\ 0       \\ 1 \end{bmatrix}, &
  \br_{\imh,j}^3 = \br_{i,j}^{3}   & = \begin{bmatrix}-1 \\ Z_{i,j} \\ 0 \end{bmatrix}
\end{align}

In the nonlinear case, the Riemann solution also consists of three waves, one of
which is a stationary shear wave with zero velocity.  The other
two waves may be rarefaction waves or shock waves, but 
the solution cannot include transonic rarefactions, since the $p$-system
is derived in a Lagrangian frame.  
For reasons of computational efficiency, we will use an approximate all-shock
solver.
Shock waves in the Riemann solution correspond to traveling discontinuities
that are proportional to $\br_{i-1,j}^1$ (left-going) or $\br_{i,j}^3$
(right-going).  Equations for an exact all-shock Riemann solution can be
derived in a straightforward manner by considering the Rankine-Hugoniot
conditions.  However, these lead to a coupled nonlinear system that is expensive
to solve numerically.  Instead, we again follow \cite{leveque2002finitePaper}
and approximate the solution further by replacing the exact wave speeds with
the sound speeds in each cell.  The waves themselves are found following
the $f$-wave approach \cite{bale2003wave}, by decomposing the jump in the normal flux
in terms of the eigenvectors \eqref{imh_eigenvectors}:
\begin{align*}
  \bF_{i,j} - \bF_{i-1,j} & = \beta_{\imh,j}^{1}\br_{\imh,j}^{1}+\beta_{\imh,j}^{2}\br_{\imh,j}^{2}+\beta_{\imh,j}^{3}\br_{\imh,j}^{3} \\
   & = \Zwave_{\imh,j}^{1}+\Zwave_{\imh,j}^{2}+\Zwave_{\imh,j}^{3},
 \end{align*}
where $\bF_{ij}=\bff(\bQ_{ij})$.
The second wave is not needed in the numerical solution, since it has
velocity zero and thus does not affect the solution value in either cell.
The other two waves are accumulated into left- and right-going fluctuations, 
which are quantities that consider the net effect of all left- and right-going waves respectively:
\begin{align}
    \amdQ_{\imh,j} & = \Zwave_{\imh,j}^1 \\
    \apdQ_{\imh,j} & = \Zwave_{\imh,j}^3.
\end{align}

The normal Riemann solver for the y-direction uses the same approach.  To solve
the Riemann problem at $(x_i,y_\jmh)$, in place of \eqref{imh_eigenvectors} we use
the eigenvectors
\begin{align} \label{jmh_eigenvectors}
  \br_{i,\jmh}^1 = \br_{i,j-1}^{1} & = \begin{bmatrix} 1 \\ 0 \\ Z_{i,j-1} \end{bmatrix}, & 
  \br_{i,\jmh}^2 = \br^{2}         & = \begin{bmatrix} 0 \\ 1       \\ 0   \end{bmatrix}, &
  \br_{i,\jmh}^3 = \br_{i,j}^{3}   & = \begin{bmatrix}-1 \\ 0 \\ Z_{i,j}   \end{bmatrix}.
\end{align}
The normal flux difference is decomposed in terms of these eigenvectors:
\begin{align*}
 \bG_{i,j} - \bG_{i,j-1}  & = \beta_{i,\jmh}^{1}\br_{i,\jmh}^{1}+\beta_{i,\jmh}^{2}\br_{i,\jmh}^{2}+\beta_{i,\jmh}^{3}\br_{i,\jmh}^{3} \\
   & = \Zwave_{i,\jmh}^{1}+\Zwave_{i,\jmh}^{2}+\Zwave_{i,\jmh}^{3},
 \end{align*}
where $\bG_{ij}=\bg(Q_{ij})$. Finally, the waves are accumulated into 
up- and down-going fluctuations:
\begin{align}
    \bmdQ_{i,\jmh} & = \Zwave_{i,\jmh}^1 \\
    \bpdQ_{i,\jmh} & = \Zwave_{i,\jmh}^3.
\end{align}
As observed in \cite{leveque2002finitePaper} for the 1D $p$-system, this
approach leads to efficient approximate Riemann solvers that are conservative
and provide good accuracy at least for weakly nonlinear problems.

  \subsection{Transverse Riemann solvers\label{sec:transverse}}

The multidimensional Clawpack algorithm makes use of a transverse Riemann
solver that decomposes the horizontally traveling waves into up- and down-going
corrections and the vertical traveling waves into right- and left-going
corrections.  These second-order corrections capture the effect of corner
transport.  

Transverse corrections to the vertical fluctuations are computed by decomposing the horizontal 
fluctuations into up- and down-going parts using the eigenvectors \eqref{jmh_eigenvectors}:
\begin{equation}
  \mathcal{A}^{\pm}\Delta\bQ_{i-\frac{1}{2},j} = \gamma_{i,j-\frac{1}{2}}^{1}\bm{r}_{i,j-\frac{1}{2}}^{1}+\gamma_{i,j-\frac{1}{2}}^{2}\bm{r}_{i,j-\frac{1}{2}}^{2}+\gamma_{i,j-\frac{1}{2}}^{3}\bm{r}_{i,j-\frac{1}{2}}^{3}.
\end{equation}
The speed of the waves in the vertical direction determines the amount
of horizontal fluctuation that must be added to the vertical one.
These speeds are given by the eigenvalues of $\mB_{i,\jmh}$, which are 
$s_{i,j-\frac{1}{2}}^{1}=-c_{i,j-\frac{1}{2}}$, $s_{i,j-\frac{1}{2}}^{2}=0$ and
$s_{i,j-\frac{1}{2}}^{3}=c_{i,j-\frac{1}{2}}$. Finally, the corrections to the
vertical fluctuations are: 
\begin{align*}
  \mathcal{B}^{-}\mathcal{A}^{\pm}\Delta\bQ_{i-\frac{1}{2},j} & = s_{i,j-\frac{1}{2}}^{1}\gamma_{i,j-\frac{1}{2}}^{1}\bm{r}_{i,j-\frac{1}{2}}^{1},\\
  \mathcal{B}^{+}\mathcal{A}^{\pm}\Delta\bQ_{i-\frac{1}{2},j} & = s_{i,j-\frac{1}{2}}^{3}\gamma_{i,j-\frac{1}{2}}^{3}\bm{r}_{i,j-\frac{1}{2}}^{3}.
\end{align*}

The transverse corrections to the horizontal fluctuations are obtained in an
analogous way. Those corrections are: 
\begin{align*}
  \mathcal{A}^{-}\mathcal{B}^{\pm}\Delta\bQ_{i,j-\frac{1}{2}} & = s_{i-\frac{1}{2},j}^{1}\gamma_{i-\frac{1}{2},j}^{1}\bm{r}_{i-\frac{1}{2},j}^{1},\\
  \mathcal{A}^{+}\mathcal{B}^{\pm}\Delta\bQ_{i,j-\frac{1}{2}} & = s_{i-\frac{1}{2},j}^{3}\gamma_{i-\frac{1}{2},j}^{3}\bm{r}_{i-\frac{1}{2},j}^{3}.
\end{align*}

  \subsection{Accuracy tests}

    For the nonlinear, variable-coefficient system studied in the present work, exact
    solutions are not available and error and convergence estimates are
    difficult to obtain because the degree of regularity of solutions is not known.
    In this section some measure of the accuracy of the numerical solutions
    is obtained by conducting self-convergence tests on small problems
    with the same qualitative features as the problems of interest.
    The principal purpose of these tests is to determine the grid resolution necessary
    to ensure small relative errors and qualitatively accurate results.
    Detailed numerical analysis of the schemes' convergence behavior for these
    problems is beyond the scope of this work.
    In all computations here and in the rest of this work, a uniform
    cartesian grid is used with $\Dx=\Dy=h$.

    \subsubsection{Sinusoidal medium}

Considering the nonlinear constitutive relation \eqref{psystem}
and the material properties $\rhoB=\KB=10$,
we perform a self-convergence study in the sinusoidal medium \eqref{sinusoidal} using
the Clawpack discretization. The initial condition is given by \eqref{eq: initial condition}.
The computational domain is restricted to the positive 
quadrant by imposing reflecting boundary conditions at the left and bottom. 
The stress at $t=3$ on a grid with $h^*=1/480$ is used as reference solution.
Table \ref{tab: Self-convergence_sinusoidal-classic unsplit} shows
self-convergence rates as well as relative errors, computed by: 
\begin{equation}\label{eq: relative error}
  E= \frac{h^*}{h}\cdot\frac{||\sigma-\sigma^*||_2}{||\sigma^*||_2},
\end{equation}
where $\sigma$ is the computed solution and $\sigma^*$ is the solution on
the fine grid.  The observed convergence rate is roughly second order.

\begin{table}[H]
  \begin{centering}
  \subtable{
    \begin{tabular}{c|c|c}
      \hline 
        $\frac{1}{h}$ & $L^2$ Error & Rate\tabularnewline
      \hline
      \hline 
        80 & 1.737x10$^{-3}$ & ---\tabularnewline
      \hline 
        120 & 8.943x10$^{-4}$ & 1.638 \tabularnewline
      \hline 
        160 & 5.306x10$^{-4}$ & 1.814 \tabularnewline
      \hline 
        240 & 2.235x10$^{-4}$ & 2.132 \tabularnewline
      \hline
    \end{tabular}}
  \par
  \end{centering}
  \caption{Self-convergence test for the sinusoidal medium using the Clawpack discretization. 
    \label{tab: Self-convergence_sinusoidal-classic unsplit}}
\end{table}

    \subsubsection{Checkerboard medium}

Next we consider self-convergence of the solution of \eqref{psystem}
in the checkerboard medium \eqref{checkerboard} using the Clawpack discretization and  SharpClaw. We use the
nonlinear constitutive relation
given by \eqref{eq: nonlinear const relation} and material
parameters $\rhoB=\KB=5$.
The initial condition is given by \eqref{eq: initial condition}. 
The stress is computed at $t=3$ and the solution on a grid 
with $h=1/480$ is taken as the reference solution
for the self-convergence study.
Tables \ref{tab: Self-convergence_checker-classic unsplit} and 
\ref{tab: Self-convergence_checker-sharpclaw} show the 
convergence rates and relative errors following \eqref{eq: relative error} 
using Clawpack and SharpClaw respectively.
The achieved order of convergence is between 1 and 2 for both discretizations,
which seems reasonable given the discontinuous coefficients.
The SharpClaw solution is noticeably more accurate, even though the two
methods exhibit similar convergence rates.

\begin{table}[H]
  \begin{centering}
   \subtable[Clawpack\label{tab: Self-convergence_checker-classic unsplit}]{
    \begin{tabular}{c|c|c}
      \hline 
        $\frac{1}{h}$ & $E$ & $p$ \\
      \hline
      \hline 
        80 & 2.935x10$^{-2}$ & --- \\
      \hline 
        120 & 1.921x10$^{-2}$ & 1.045 \\
      \hline 
        160 & 1.358x10$^{-2}$ & 1.205 \\
      \hline 
        240 & 7.293x10$^{-3}$ & 1.533 \\
      \hline
    \end{tabular}}
   \subtable[SharpClaw\label{tab: Self-convergence_checker-sharpclaw}]{
    \begin{tabular}{c|c|c}
      \hline 
        $\frac{1}{h}$ & $E$ & $p$ \\
      \hline
      \hline 
        80 & 1.223x10$^{-2}$ & --- \\
      \hline 
        120 & 7.463x10$^{-3}$ & 1.218 \\
      \hline 
        160 & 5.035x10$^{-3}$ & 1.367 \\
      \hline 
        240 & 2.561x10$^{-3}$ & 1.667 \\
      \hline
    \end{tabular}}
  \par
  \end{centering}
  \caption{Self-convergence test for the checkerboard medium using 
  (a) Clawpack and (b) SharpClaw.
  \label{tab: Self-convergence_checker}}
 \end{table}

\section{Computational results\label{sec:results}}

  \subsection{Formation of solitary waves}
We now investigate in detail the solitary wave trains mentioned in Section \ref{sec:introduction}.
These waves are found to arise in both the checkerboard and the sinusoidal domain.  We
take a symmetric Gaussian hump \eqref{eq: initial condition} as initial stress
and solve the $p$-system \eqref{psystem-2nd-order} with the nonlinear constitutive relation
\eqref{eq: nonlinear const relation}.
Reflecting boundary conditions are used at the left and bottom boundaries in order to
take advantage of symmetry and compute only in the first quadrant.  The mesh is
uniform with $h=\frac{1}{240}$.  Based on the accuracy studies above and additional
tests, we expect this to be sufficient computational resolution to obtain qualitatively
accurate solutions.

\begin{figure}
 \begin{centering}
  \subfigure[Checkerboard medium, $t=90$, $h=1/240$\label{Flo: Stress using checker at t=90}]{
        \includegraphics[scale=0.3]{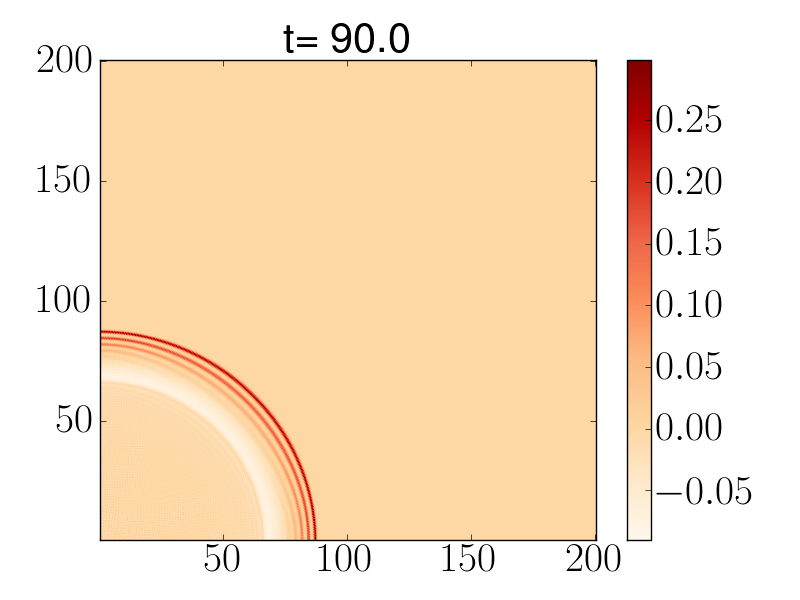}
        \includegraphics[scale=0.3]{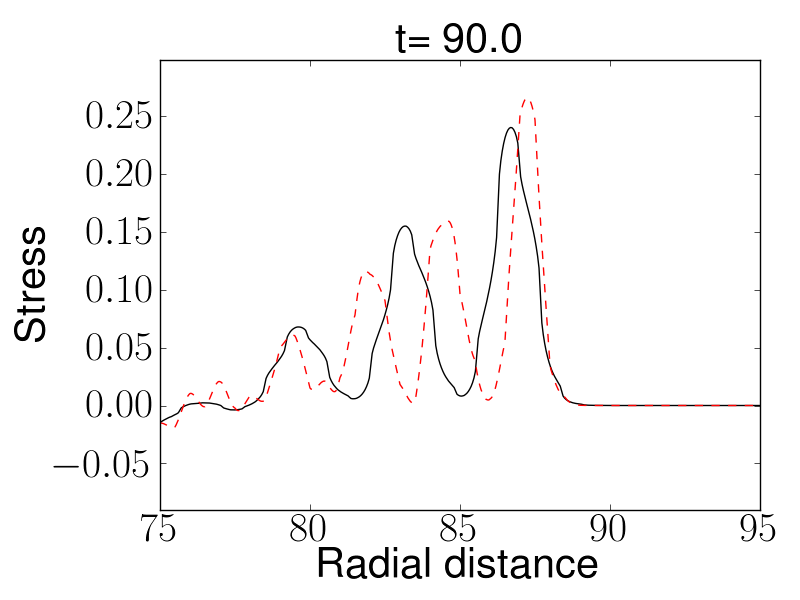}}
  \subfigure[Checkerboard medium, $t=200$, $h=1/240$\label{Flo: Stress using checker at t=200}]{
        \includegraphics[scale=0.3]{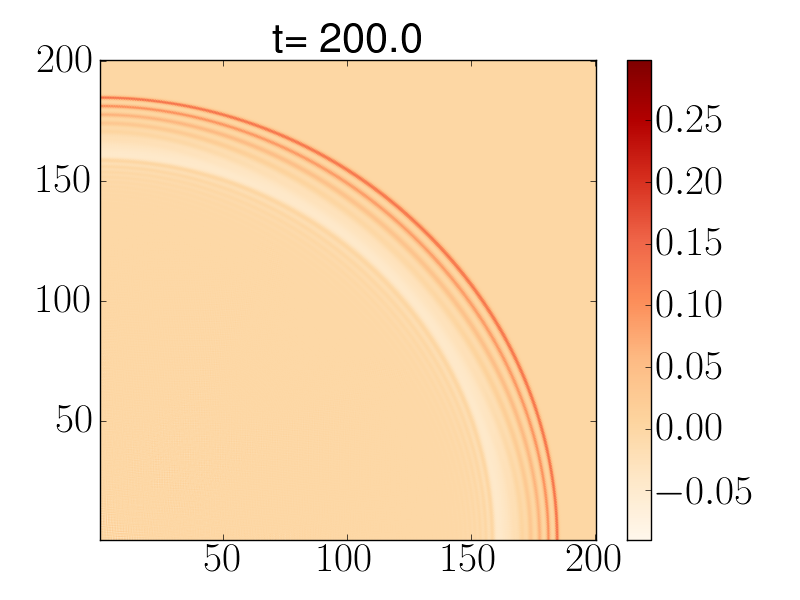}
        \includegraphics[scale=0.3]{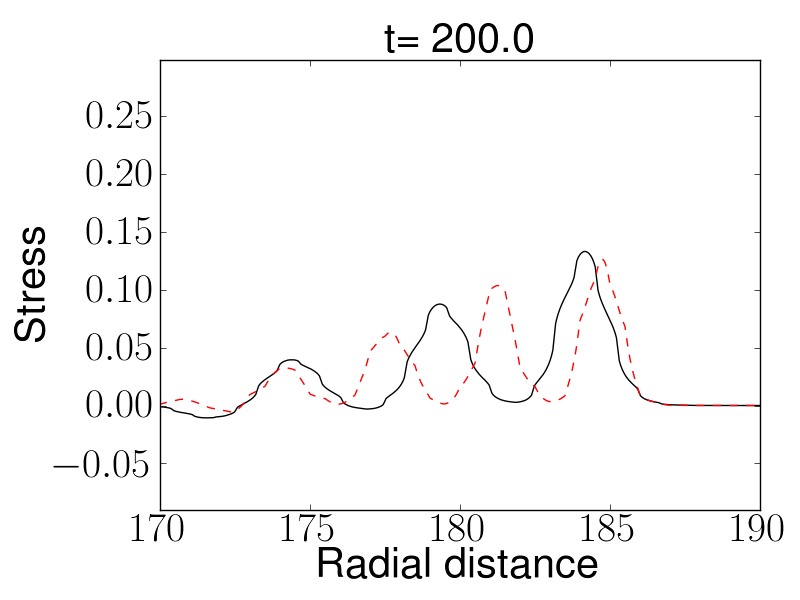}}      
  \subfigure[Sinusoidal medium, $t=90$, $h=1/240$\label{Flo: Stress using sinusoidal at t=90}]{
        \includegraphics[scale=0.3]{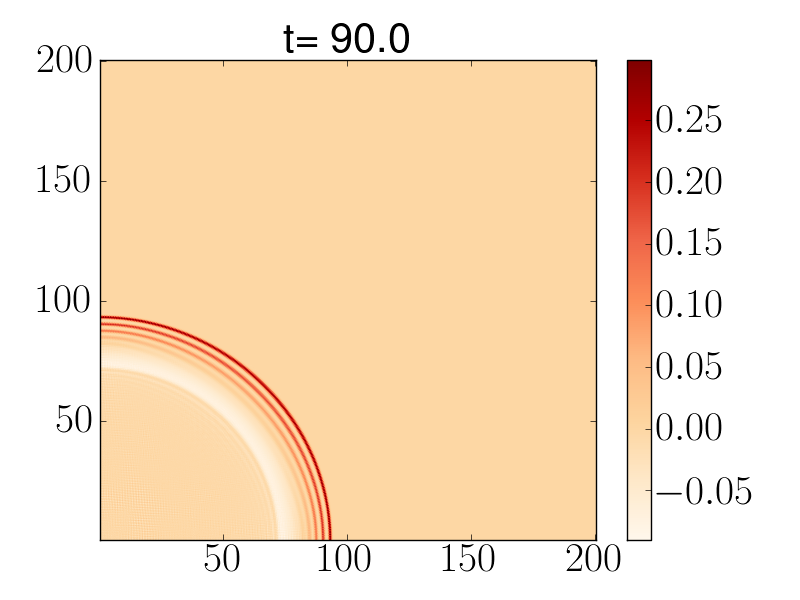}
        \includegraphics[scale=0.3]{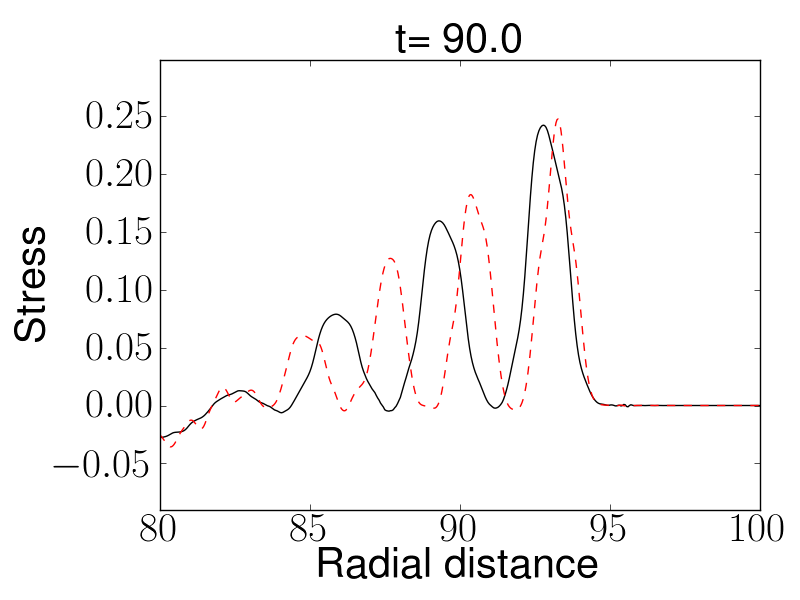}}
  \subfigure[Sinusoidal medium, $t=200$, $h=1/120$\label{Flo: Stress using sinusoidal_t200}]{
        \includegraphics[scale=0.3]{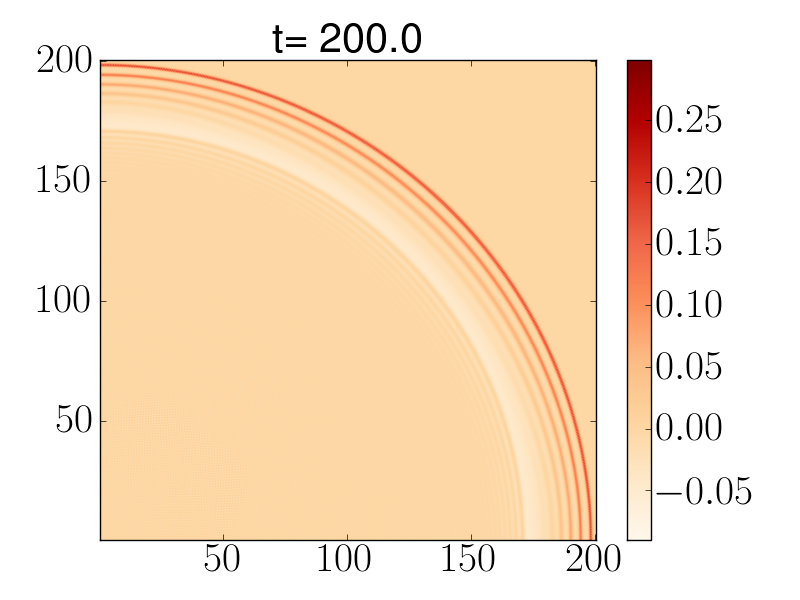}
        \includegraphics[scale=0.3]{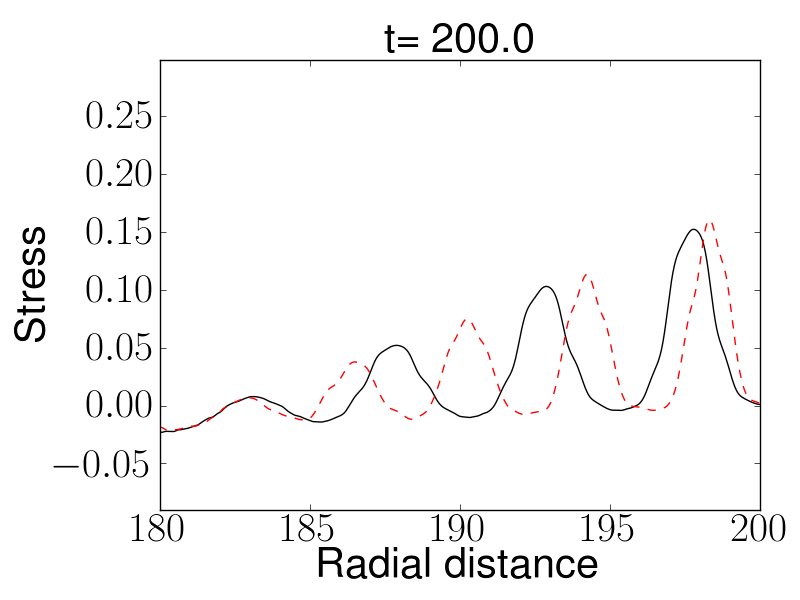}}
 \caption{Stress at $t=90$ and $t=200$ for solitary wave trains in the sinusoidal and checkerboard media.  
           The slice plots (on the right) show results along the lines $y=x$ (solid black line) and $y=0$ (dashed red line).\label{fig:wave-formation}}
 \end{centering}
\end{figure}

For the checkerboard domain, we take $\rhoB =\KB= 5$ and compute the solution
using SharpClaw with $h=\frac{1}{240}$. For the sinusoidal domain, we take $\rhoB =\KB= 10$ and compute
the solution using Clawpack with $h=\frac{1}{120}$.
Figure \ref{fig:wave-formation} shows the stress at
$t=90$ and $t=200$ for both materials, including slices along the lines $y=x$ and $y=0$. 
The persistence of the solitary waves after long times strongly suggests
that they are stable, attracting solutions.

Note that the solitary wave front is noticeably non-circular and
that the cross-sectional shape of the solitary wave pulse varies with respect
to the angle of propagation.

  \subsection{Entropy evolution and shock formation}
As discussed in the introduction,
solutions of nonlinear hyperbolic PDEs generically develop shock discontinuities,
leading to irreversibility and entropy decay.  It has been observed that spatially
varying materials can inhibit the formation of shocks \cite{fouque2004shock,Ketcheson_LeVeque_2011}.
Meanwhile, dispersive nonlinear wave equations commonly exhibit solutions that are regular for
all times. Here we are solving a PDE without dispersive terms; nevertheless,
the spatially varying coefficients create reflections,
which yield effective dispersion.

\begin{figure}
 \begin{centering}
  \subfigure{\includegraphics[scale=0.3]{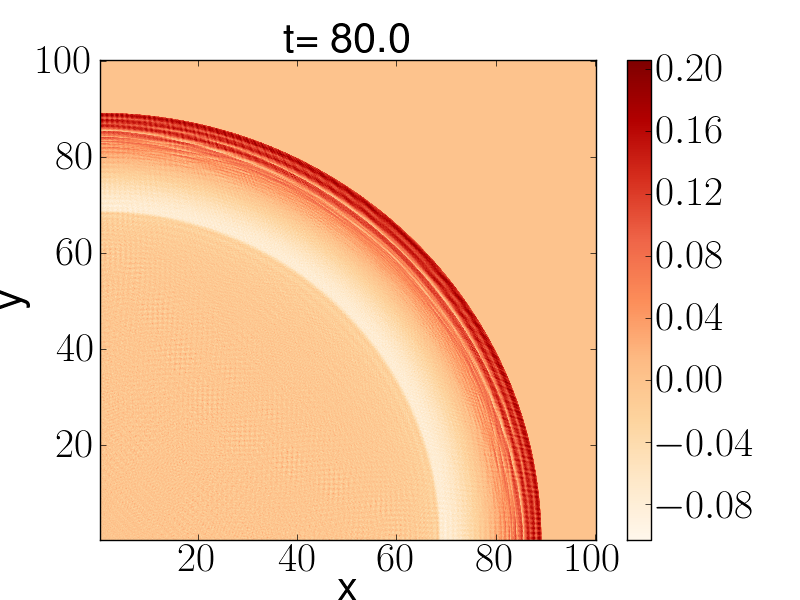}\label{Flo: introduction psystem_het nonlinear small ratio}}
  \subfigure{\includegraphics[scale=0.3]{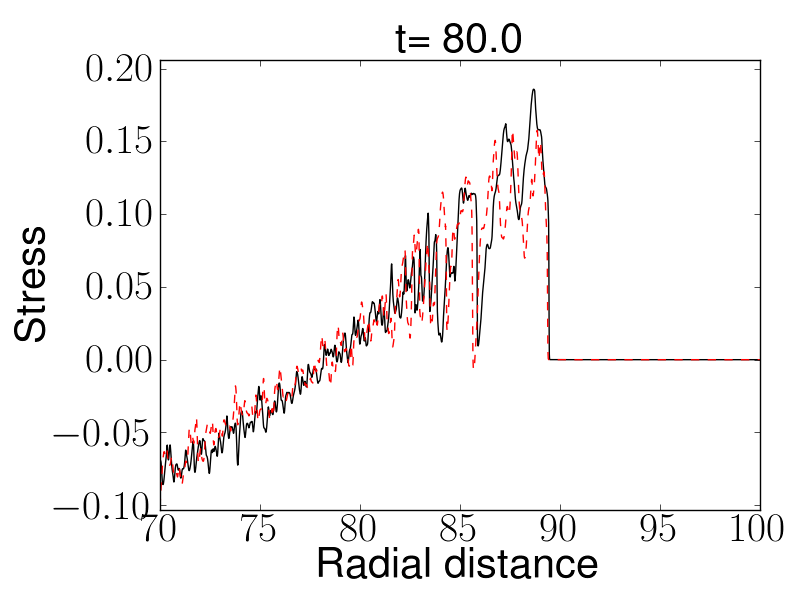}\label{Flo: introduction psystem_het nonlinear small ratio slices}}
  \par
 \end{centering}
 \caption{(a) Stress at $t=90$ in a checkerboard medium with small impedance ratio: $\rhoB =\KB= 1.5$. \label{Flo: small vs large impedance ratio} (b) Slice plots showing results along the lines $y=x$ (solid black line) and $y=0$ (dashed red line).}
\end{figure}

No visible shocks appear in the solutions shown in Figure \ref{fig:wave-formation}.
On the other hand, figure \ref{Flo: small vs large impedance ratio} shows that
if the impedance ratio
is small enough, the solution develops shocks.  This is in qualitative agreement with the theory
for 1D systems developed in \cite{Ketcheson_LeVeque_2011}.
We now study shock formation in the 2D stegotons by considering the
entropy evolution and the reversibility of the solution. 

An entropy function for a hyperbolic PDE is a function that is conserved
while the solution is smooth, but decreases in time when shocks form \cite{leveque2002finite}.
Thus, entropy decay is a signature of shock formation.  A suitable entropy function 
for the 2D p-system \eqref{psystem-2nd-order} is the total energy:
\begin{equation}
 \eta(u,\epsilon)=\int_{\bm{x}}\left( \frac{1}{2}\rho(\bm{x})u^2+\int_0^\epsilon\sigma(s,\bm{x})ds\right) dx.
\end{equation}

We now compare the solution in a homogeneous domain 
with $\rho=K=1$ (shown in Figure \ref{Flo: Stress in homogeneous and
heterogeneous nonlinear}) to that obtained in
the sinusoidal domain with $\rhoB=\KB=10$ (shown in Figure 
\ref{Flo: Stress using sinusoidal at t=90}).
The difference in the behavior is clear: the heterogeneity introduces
reflections that effectively yield dispersion, breaking the initial profile
into solitary waves. For the homogeneous domain, no reflections are present;
consequently, there is no dispersion to prevent shock formation. In Figure
\ref{Flo: Entropy hom vs het}, we see the entropy evolution for both cases.
 For the simulation with
homogeneous domain, the entropy starts to decrease as soon as the shock starts
to happen. On the other hand, the entropy using the sinusoidal type medium
remains almost constant. In both cases we use a resolution given by
$h=\frac{1}{240}$. 

\begin{figure}
 \begin{centering}
  \subfigure{\includegraphics[scale=0.3]{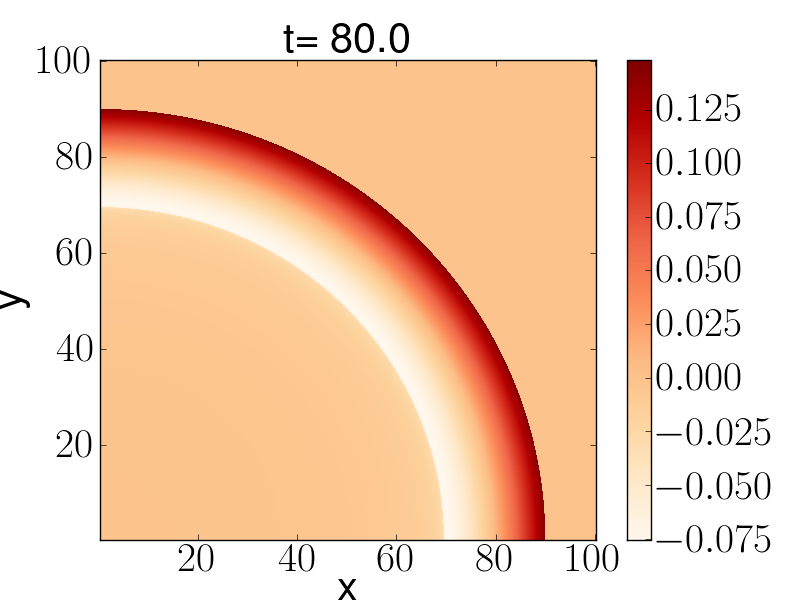}\label{Flo: Stress in homogeneous nonlinear}}
  \subfigure{\includegraphics[scale=0.3]{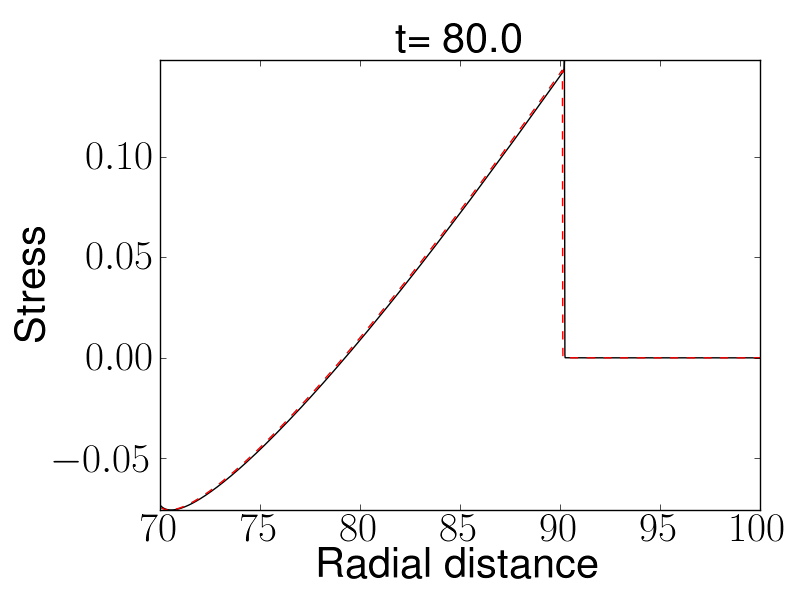}\label{Flo: Stress in heterogeneous nonlinear}}
  \par
 \end{centering}
 \caption{(a) Stress at $t=80$ in a homogeneous nonlinear medium. (b) Slice plots showing results along the lines $y=x$ (solid black line) and $y=0$ (dashed red line). \label{Flo: Stress in homogeneous and heterogeneous nonlinear}}
\end{figure}

\begin{figure}
 \begin{centering}
   \includegraphics[scale=0.3]{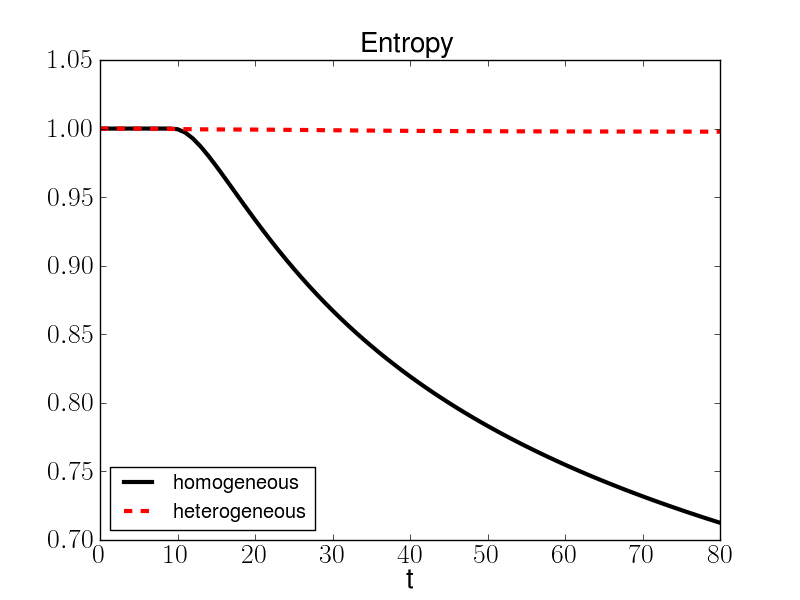}
   \par
 \end{centering}
 \caption{Entropy evolution considering a homogeneous nonlinear medium
 (solid line) and a heterogeneous nonlinear medium (dashed line). The entropy
 has been normalized relative to the initial entropy.\label{Flo: Entropy
 hom vs het}}
\end{figure}

\begin{figure}
 \begin{centering}
  \subfigure[Entropy evolution up to $t=80$]{\includegraphics[scale=0.3]{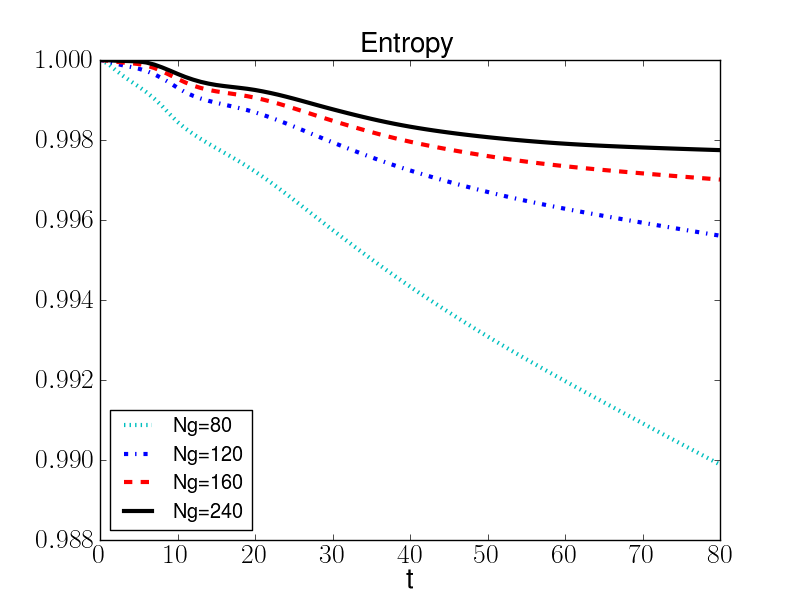}\label{Flo: Entropy upto t80}}
  \subfigure[Entropy evolution up to $t=10$]{\includegraphics[scale=0.3]{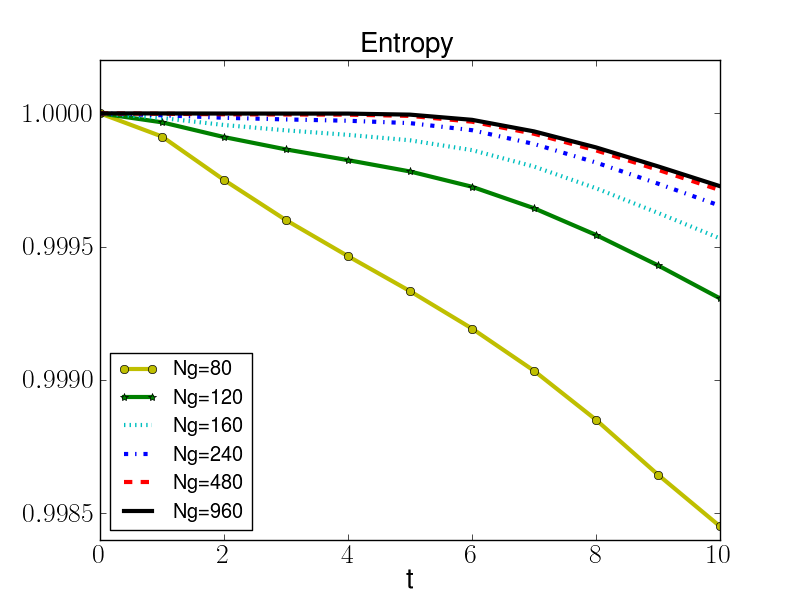}\label{Flo: Entropy upto t10}}
  \par
 \end{centering}
 \caption{Entropy evolution up to (a) $t=80$ and (b) $t=10$ considering (a) $\frac{1}{h}=80, 120, 160, 240$ and (b) $\frac{1}{h}=80, 120, 160, 240, 480, 960$. The entropy
 has been normalized relative to the initial entropy.}
\end{figure}

Indeed, in both cases, there is loss of entropy due to numerical dissipation, but this converges to zero, as the grid is refined. To better determine if there is shock formation in the sinusoidal case, we present, in figure \ref{Flo: Entropy upto t80}, the normalized entropy evolution considering different resolutions given by $\frac{1}{h}=80, 120, 160, 240$. We see the entropy lost is indeed small. However, it is evident that the rate of entropy lost changes before $t=10$ making the entropy lost not to converge to zero as the grid is refined. This suggests the existence of a shock. To corroborate this, we focus up to $t=10$, study the entropy for even finer grids and see its convergence. In figure \ref{Flo: Entropy upto t10}, we show the entropy evolution up to $t=10$ for different resolutions given by $\frac{1}{h}=80, 120, 160, 240, 480, 960$. We see the entropy at $t=10$ converges to a value different than the initial entropy; i.e., there is a loss of entropy, which corroborates there existence of a shock. 

\subsection{Interaction of cylindrical stegotons}
\begin{figure}
 \begin{centering}
  \subfigure[Spatial distribution of stress.]{
  \begin{tabular}{cc}
   \includegraphics[scale=0.35]{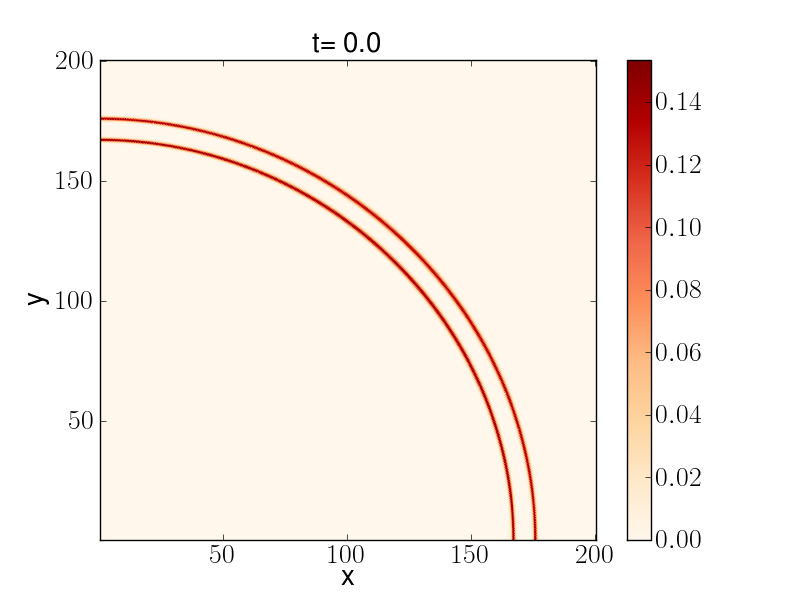} &
   \includegraphics[scale=0.35]{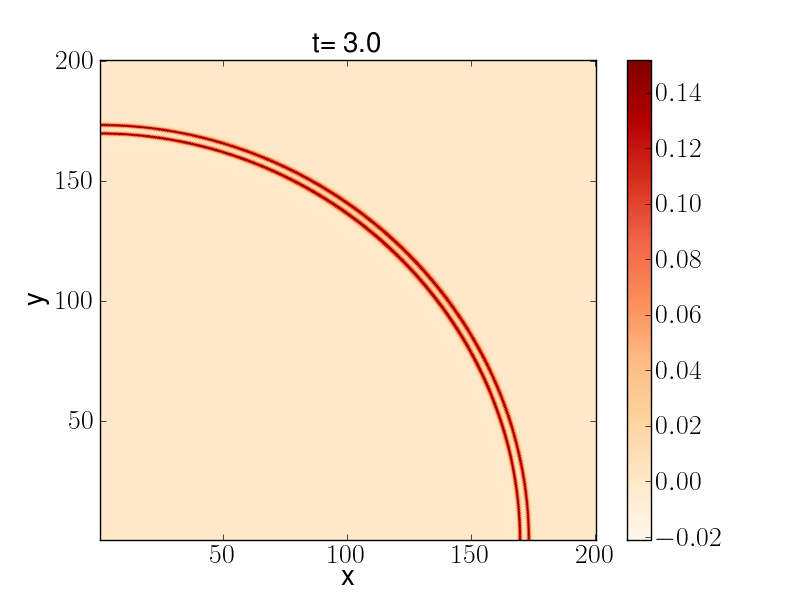} \\
   \includegraphics[scale=0.35]{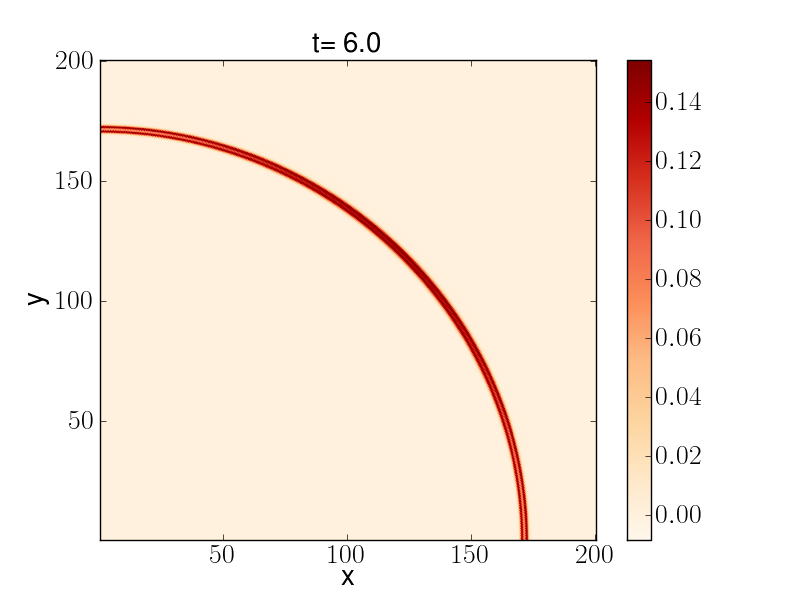} &
   \includegraphics[scale=0.35]{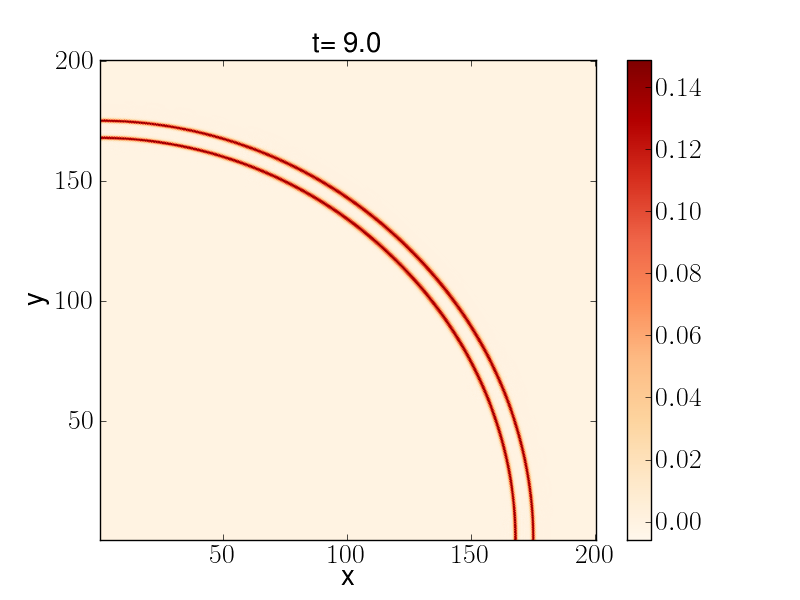} 
   \end{tabular}
  \label{fig: 2D collision_surface plot}}
 \subfigure[Slices at $y=x$.]{
 \begin{tabular}{cc}
  \includegraphics[scale=0.35]{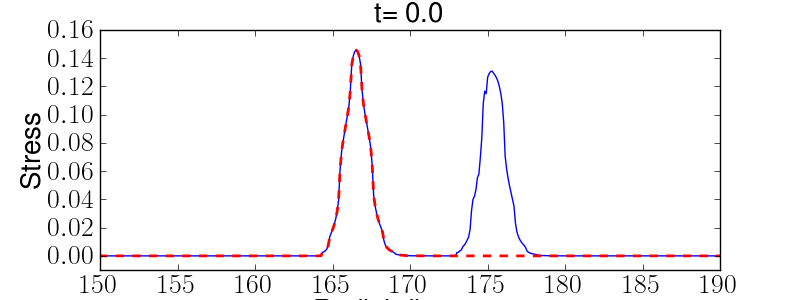} &
  \includegraphics[scale=0.35]{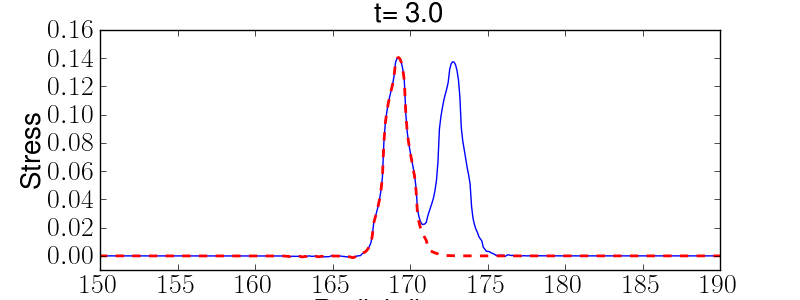} \\
  \includegraphics[scale=0.35]{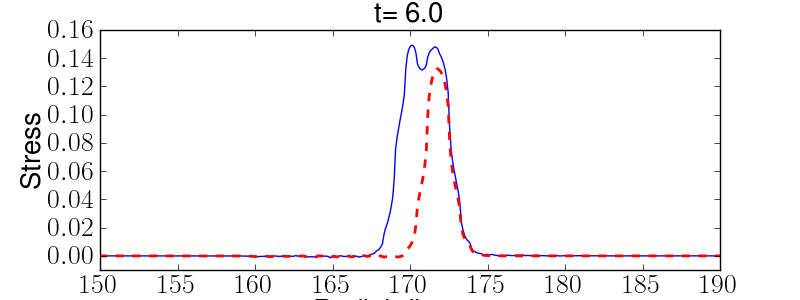} &
  \includegraphics[scale=0.35]{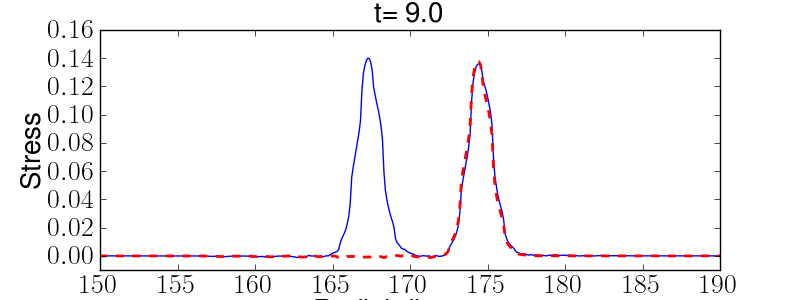}
  \end{tabular}
   \label{fig: 2D collision_slices_45}}
 \end{centering}
 \caption{2D collision of cylindrical stegotons for $t=0, 3, 6, 9$. The
            stegotons are initially moving toward each other.\label{fig: 2D collision}}
\end{figure}

In order to study the interaction of these waves, we start with the solution
depicted in Figure \ref{Flo: Stress using checker at t=200}.  We attempt to
isolate the leading pulse in that solution by setting the solution to zero
outside of a narrow band. Let $\bq$ denote the solution of \eqref{psystem} 
and $\bqh$ denote the corresponding isolated
solution. Then $\bqh$ is given by: 

\begin{align}
\bqh(x_i,y_j) & = \begin{cases}
\bq(x_i,y_j) & \text{ if } x_i \ge X(j), \\
0        & \text{ if } x_i < X(j).
\end{cases}
\end{align}
Here $X(j)$ is chosen so as to separate the leading stegoton from
everything to the left as well as possible.  Specifically, for 
each $j$, $X(j)$ is the right-most local minimum of $\bq$:
\begin{align}
X(j) & = \max_i \{x_i : S_{ij} < S_{i+1,j}; S_{ij} < S_{i-1,j} \},
\end{align}

where $S$ is the discretization of the stress $\sigma$. 
This isolation
is imperfect since it seems clear that the tails of the solitary waves still
overlap.  Obtaining completely separated solitary waves would require even greater
computational resources or a different modeling approach.  By extracting just
the leading pulse in this manner at two different times ($t=180$ and $t=190$),
we obtain a pair of nearly isolated solitary waves.  The interaction simulation
is initialized using the sum of these two solutions as initial condition, but
with the velocity negated in the $t=190$ solution so that it will propagate inward.  
Figure \ref{fig: 2D collision_surface plot} shows surface plots of the solution at 
$t=0, 3, 6, 9$. Corresponding slices at $y=x$ are shown in Figure 
\ref{fig: 2D collision_slices_45} (solid line). For comparison, we also
simulate just the outer pulse (without the inner pulse) over the same time
(dotted line).

\begin{figure}
 \begin{centering}
  \subfigure[1D collision with same direction.]{
   \includegraphics[scale=0.35]{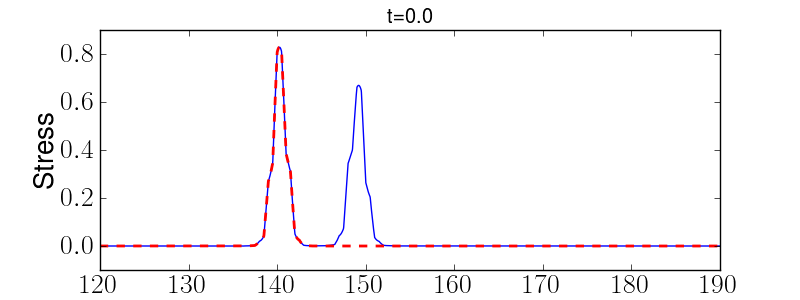}
   \includegraphics[scale=0.35]{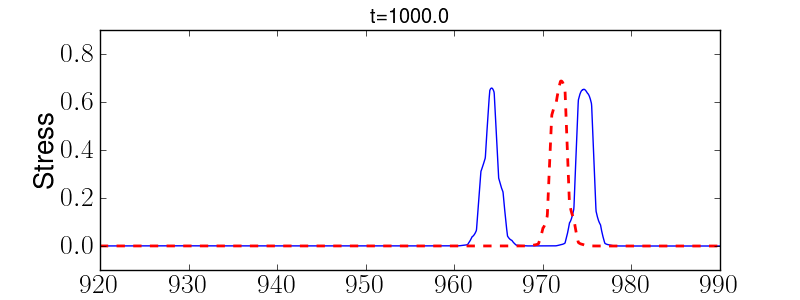}
    \label{fig: 1D collision_same direction}}
 \subfigure[1D collision with opposite direction.]{
  \includegraphics[scale=0.35]{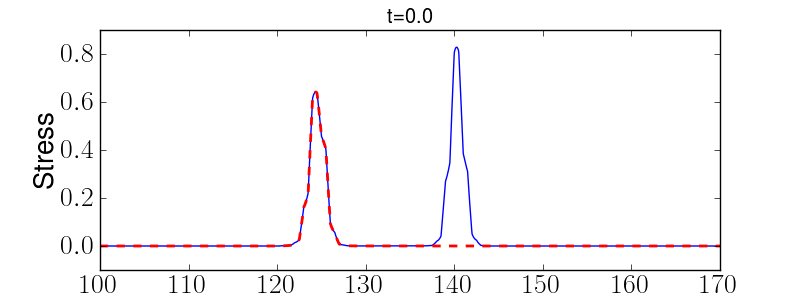}
  \includegraphics[scale=0.35]{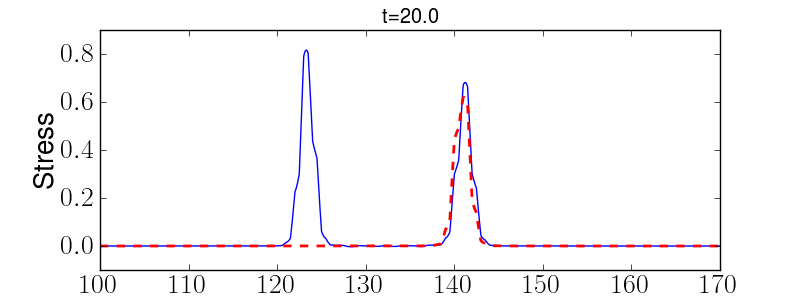}
   \label{fig: 1D collision_opposite direction}}
  \par
 \end{centering}
 \caption{Interaction of 1D stegotons. The dashed line shows the stegoton originally at the left propagating on its own.\label{fig: 1D interaction}}
\end{figure}

The results are typical of solitary waves; the two pulses retain their
identity after the interaction.  Furthermore, the position of the outward-going
pulse after the interaction is essentially the same as in the simulation
with no interaction.  The apparent lack of a noticeable phase shift seems
to be a result of the very short interaction time for this head-on collision.
The same effect can be seen in one-dimensional stegoton interactions, as shown
in Figure \ref{fig: 1D interaction}.  Two stegotons traveling in the same 
direction exhibit a phase shift after interaction, but no phase shift is 
discernible after a head-on interaction.  In order to simulate the interaction
of two cylindrical stegotons traveling in the same direction (and the 
associated phase shift), much greater computational resources would be required.

Some small oscillations are visible in the final solution in Figure \ref{fig: 2D collision}.
Since similar oscillations are seen in the solution obtained with just the
outward-going pulse, it seems that these may be attributed to the fact that the
tails of the solitary waves were not accurately captured in the initial condition.
However, it might also be that some small oscillations are radiated during the
interaction.

\section{Conclusion}
We have performed direct numerical simulations of the multiscale behavior
of nonlinear waves in non-dispersive periodic materials.  The largest simulations
involve $6.9\times 10^9$ unknowns and were performed on 16K cores of a BlueGene/P
system.
Our computational results indicate that solitary waves may arise in solutions
of two-dimensional first-order nonlinear hyperbolic PDEs with spatially
periodic coefficients, including both smooth and piecewise-constant media.
As in 1D, these waves apparently result from the combination of effective (material)
dispersion and nonlinear steepening.  The effective dispersion is a macroscopic
effect caused by reflections in variable-impedance media.  In media
with strongly varying impedance, shock formation is strongly suppressed
relative to that occurring in homogeneous media, as demonstrated by near-conservation
of entropy.  The cylindrical solitary waves that we have studied appear to
behave approximately like solitons in their interactions.

Studying the properties of these waves in detail is difficult due to
the multiscale nature of the problem.  Future work may include simulations
using massively parallel adaptive mesh refinement or more computational
resources in order to more fully isolate the solitary waves and study
longer interactions.  A complementary tool in understanding these waves
is the application of homogenization theory to multidimensional nonlinear
hyperbolic systems.  
We have focused on localized perturbations in media with uniform sound
speed and strongly varying impedance.
Many other geometries, for both the medium and the perturbation, could also be
considered.  All of these topics are the subject of ongoing work.

\bibliographystyle{plain}
\bibliography{bibs}

\end{document}